\documentclass[11pt]{article}

\usepackage{tikz}
\usepackage{subfigure}
\usepackage[english]{babel}
\usepackage{graphicx}
\usetikzlibrary{arrows, graphs}
\usepackage{amsfonts,amssymb,amsmath,latexsym,amsthm}
\usepackage{multirow}
\usepackage{bbm}
\usepackage[noabbrev,capitalize,nameinlink]{cleveref}

\newcommand \HH{\mathcal H}
\newcommand \kd{\operatorname{def}}
\newcommand\B{\text{B-}}

\newtheorem{thm}{Theorem}[section]

\newtheorem{lem}[thm]{Lemma}
\newtheorem{prop}[thm]{Proposition}
\newtheorem{prob}[thm]{Problem}
\newtheorem{conj}[thm]{Conjecture}
\newtheorem{quest}[thm]{Question}

\newtheorem{claim}{Claim}

\DeclareFontFamily{U}{mathx}{}
\DeclareFontShape{U}{mathx}{m}{n}{<-> mathx10}{}
\DeclareSymbolFont{mathx}{U}{mathx}{m}{n}
\DeclareMathAccent{\widehat}{0}{mathx}{"70}
\DeclareMathAccent{\widecheck}{0}{mathx}{"71}
\newcommand \fmz{\widecheck}
\newcommand \diam{\operatorname{diam}}
\newcommand \dist{\operatorname{dist}}
\begin{document}

\title{A Ramsey-type theorem on deficiency}
\author{Jin Sun$^a$, \quad Xinmin Hou$^{a,b,c}$\footnote{\text{e-mail: jinsun@mail.ustc.edu.cn (J. Sun), xmhou@ustc.edu.cn (X. Hou)}}\\
\small $^{a}$ School of Mathematical Sciences\\
\small University of Science and Technology of China, Hefei, Anhui 230026, China.\\
\small  $^{b}$ CAS Key Laboratory of Wu Wen-Tsun Mathematics\\
\small University of Science and Technology of China, Hefei, Anhui 230026, China.\\
\small$^c$ Hefei National Laboratory,\\
\small University of Science and Technology of China, Hefei 230088, Anhui, China
}
\date{}

\maketitle

\begin{abstract}
Ramsey's Theorem states that a graph $G$ has bounded order if and only if $G$ contains no  complete graph $K_n$ or empty graph $E_n$ as its induced subgraph. 
The Gy\'arf\'as-Sumner conjecture says that a graph $G$ has bounded chromatic number if and only if it contains no induced subgraph isomorphic to  $K_n$ or a tree $T$.
The deficiency of a graph is the number of vertices that cannot be covered by a maximum matching. In this paper, we prove a Ramsey type theorem  for deficiency, i.e., we characterize all the forbidden induced subgraphs for graphs $G$ with bounded deficiency. As an application, we answer a question proposed by Fujita, Kawarabayashi, Lucchesi, Ota, Plummer and Saito (JCTB, 2006).
	
\end{abstract}

\section{Introduction}
In this paper, all the graphs we consider are finite and simple. For positive integer $n$, write $[n]$ for $\{1,2,\dots ,n\}$.   
Let $G$ be a graph with vertex set $V(G)$ and edge set $E(G)$. 
Write $|G|=|E(G)|$ for the {\em size} of $G$ and call $|V(G)|$ the {\em order} of $G$.
For a vertex $x\in V(G)$, let $N_G(x)=\{y\in V(G):xy\in E(G)\}$ be the {\it neighborhood} of $x$. Write $d_G(x)=|N_G(x)|$ for the {\em degree} of $G$.
A {\em leaf} is a vertex of degree one in a graph. 
If there is no confusion from the context, we  omit the index $G$. 
Write $\delta(G)$ and $\Delta(G)$ for the {\em minimum degree} and {\em maximum degree} of $G$, respectively.
A graph $G$ is called {\em $k$-regular} if $\delta(G)=\Delta(G)=k$. 
For two vertices $x,y$ in a connected graph $G$, the {\em distance} between $x$ and $y$ is $\dist_G(x,y)=\min\{|P| : P \text{ is an $x$-$y$ path in $G$}\}$, and the {\it diameter} of $G$ is  
$$
\diam(G)= \max_{x,y\in V(G)}\dist(x,y).    
$$
For a subset $X\subset V(G)$, let $G[X]$ denote the subgraph of $G$ induced by $X$, and set 
$
N_G(X)=\bigcup_{x\in X}N_G(x)-X.
$ 
As usual, let $P_n,\, C_n,\, K_n$ and $E_n$ denote a {\em path}, a {\em cycle}, a {\em complete graph},  and an {\em empty graph} of order $n$, respectively.
Let $K_{s,t}$ be the {\em complete bipartite graph} with partitions of orders $s$ and $t$.  
When $s=1$, we call  $K_{1,t}$ a {\em star}.  
The {\em center} of a graph $G$ is a vertex $c\in V(G)$ with $d(c,x)\le \lceil\frac{\diam(G)}2\rceil$ for any vertex $x\in V(G)$. 
Clearly, a star $K_{1,t}$ has a unique center, and a path $P_n$ has a unique center if $n$ is odd and two centers if $n$ is even. 
One can refer to \cite{d} for any other notation not introduced here.

For graphs $H_1$ and $H_2$, we say $H_1\prec H_2$ if $H_2$ contains an induced subgraph isomorphic to $H_1$, otherwise, write $H_1\nprec H_2$. For a family $\mathcal{H}$ of graphs, we say $G$ is {\it $\mathcal{H}$-free} if  $H\nprec G$ for any $H\in \mathcal{H}$. 
For two families $\mathcal{H}_{1}$ and $ \mathcal{H}_2$, we say $\mathcal{H}_1\le \mathcal{H}_2 $ if for every $H_2\in \mathcal{H}_2$, there exists $H_1\in \mathcal{H}_1$ such that $ H_1\prec H_2$. 
The following result is a simple observation; however,  it is important.
\begin{prop}
The relation `$\le$' is transitive but not antisymmetric. Therefore, if $\mathcal{H}_1\le \mathcal{H}_2$, then every $\mathcal{H}_1$-free graph is also $\mathcal{H}_2$-free. 
\end{prop}

Using this language, 
the classical Ramsey's Theorem~\cite{r} and  its connected version can be rephrased as follows.

\begin{thm}\label{THM: Ramsey}
(A) (Ramsey's Theorem~\cite{r}, 1929) For a family $\HH$ of graphs, there is a constant $c=c(\HH)$ such that $|V(G)|< c$ for every $\HH$-free graph $G$ if and only if $\HH \le \{K_n,E_n\}$ for some positive integer $n$.

(B) (Proposition 9.4.1 in \cite{d}) For a family $\HH$ of graphs, there is a constant $c=c(\HH)$ such that $|V(G)|< c$ for every connected $\HH$-free graph $G$ if and only if $\HH \le \{K_n, K_{1,n}, P_n\}$ for some positive integer $n$. 
\end{thm}
The least integer $c=c(\HH)$ when $\HH=\{K_m, E_n\}$ in (A) of Theorem~\ref{THM: Ramsey} is denoted by $R(m,n)$, called the {\em Ramsey number}. 

In this article, we consider a more general Ramsey-type problem  as follows. Let $\mathcal{G}_c$ be the family of connected graphs. 
Given a graph parameter $\mu$ and a constant $b>0$, define  
\[
b\text{-} \B \mu=\{\HH \subset \mathcal{G}_c :  \mu(G)< b \text{ for any $\HH$-free graph $G\in\mathcal{G}_c$} \},
\]
and  
\[
\B\mu=\{\HH\subset \mathcal{G}_c : \text{there is a constant } b>0 \mbox{ such that } \HH\in b\text{-} \B \mu\}. 
\]

\begin{prob}\label{PROB: mu}
Given a graph parameter $\mu$, determine $\HH\in \B\mu$.	
\end{prob}

Let $\mathcal{F}$ be a family of graphs. An {\em (induced) $\mathcal{F}$-cover} (or {\em (induced) $\mathcal{F}$-partition}) number  of a graph $G$ is the minimum $n$ such that there exist $G_1,G_2,\cdots ,G_n$ satisfying every vertex of $G$ is covered by at least (or exactly) one $G_i$, where each $G_i$ is a (induced) subgraph of $G$ isomorphic to a member of $\mathcal{F}$.
Write $\mathcal{F}\text{-c}(G)$ (or $\mathcal{F}\text{-p}(G)$) and $\text{ind-}\mathcal{F}\text{-c}(G)$ (or $\text{ind-}\mathcal{F}\text{-p}(G)$) for the $\mathcal{F}$-cover (or $\mathcal{F}$-partition) number and induced $\mathcal{F}$-cover (or induced $\mathcal{F}$-partition) number  of a graph $G$, respectively.
Different graph parameters can be formulated using this language.
For example, the {\em chromatic number} $\chi(G)=\text{ind-}\mathcal{F}\text{-p}(G)$, where $\mathcal{F}$ consists of empty graphs;
the {\em order} $\text{ord}(G)=\mathcal{F}\text{-p}(G)$, where $\mathcal{F}$ consists of $K_1$; 
and the domination number $\gamma(G)=\mathcal{F}\text{-p}(G)$, where $\mathcal{F}$ consists of stars. 
The Gy\'arf\'as-Sumner conjecture \cite{g} can be restated as follows. This conjecture is still open in general although it has been widely considered in literatures. One can refer to a survey \cite{ss}  about the progress of this conjecture.
\begin{conj}[Gy\'arf\'as-Sumner conjecture \cite{g}]
	For a finite family $\HH$ of connected graphs, $\HH\in \B \chi$ if and only if $\HH \le \{K_n, T\}$, where  $T$ is a tree and $n$ is a positive integer.
\end{conj}

In the following, we list more results on different parameters of \cref{PROB: mu}. 
\begin{description}
\item[(1)]	The connected version of the Ramsey theorem ((B) in Theorem~\ref{THM: Ramsey}) can be rephrased as $\HH\in\B$ord if and only if $\HH\le\{K_n, K_{1,n}, P_n\}$.
\item[(2)] Chiba and Furuya  determined $\HH\in\B\mu$ when $\mu$ is the $\mathcal{F}$-cover/partition number for $\mathcal{F}$ being the set of paths in \cite{cf22}, 
and when $\mu$ is the induced $\mathcal{F}$-cover/partition number for $\mathcal{F}$ being the set of paths and stars in \cite{cf23}.
\item[(3)] Choi, Furuya, Kim and Park \cite{cfkp} determined $\HH\in\B\mu$ when $\mu$ is the induced matching number (the maximum size  of an induced 1-regular subgraph of a graph) or matching number (the maximum size  of a 1-regular subgraph of a graph).

\item[(4)] Furuya \cite{f}  determined $\HH\in\B\gamma$, where $\gamma$ is the domination number of a graph.
\item[(5)] Lozin \cite{l17} determined $\HH\in\B\mu$, where $\mu$ is the neighborhood diversity or VC-dimension (two parameters related to the property of the closed neighborhood of vertices, one can refer to \cite{l17} for details).
\item[(5)] Lozin and Razgon \cite{lr} determined $\HH\in\B\mu$, where $\mu(G)$ is the tree-width of graph $G$, which is the smallest integer such that graph $G$ has a tree-decomposition into parts with at most $k+1$ vertices.
\item[(6)] Kierstead and Penrice \cite{kp}, and Scott, Seymour, and Spirkl \cite{sss} determined $\HH\in\B\partial$, where $\partial(G)$, the degeneracy of graph $G$, is the smallest integer $d$ such that every nonnull subgraph of $G$ has a vertex of degree at most $d$.   
\end{description}

In this paper, we focus on the deficiency of graphs.  
Let $M$ be a matching of a graph $G$. 
For a vertex $v\in V(G)$, we call $M$ saturates $v$ if $v\in V(M)$; otherwise, $M$ misses $v$. The {\em matching number} $\nu(G)$ of graph $G$ is the cardinality of a maximum matching of $G$.  The {\it deficiency} of graph $G$, denoted by $\kd(G)$, is defined as 
$$\kd(G)=|V(G)|-2\nu(G).
$$
According to the definition, $\kd(G)\equiv |V(G)| \pmod 2 $. Note that $\kd(G)=\mathcal{F}$-$p(G)-1$ where $\mathcal{F}=\{K_1, mK_2: m\ge 1\}$.
A vertex set $X$ in a graph $G$ is called an {\em independent set}, or {\em stable set} if $G[X]$ is an empty graph. 
For a graph $G$, let $\alpha(G)$ be the independence number of $G$. We call a vertex $v\in V(G)$ the {\em claw center} if $\alpha(G[N(v)])\ge 3$. We call $G$ {\em claw-free} if $G$ contains no claw center.
The classical results of (near-)perfect matching on claw-free graph due to Las Vergnas \cite{l}, Sumner \cite{s},  J\"unger, Pulleyblank and Reinelt~\cite{j} can be described as follows.
 
 \begin{thm}[\cite{j,l,s}]\label{claw-free}
Every connected claw-free graph has deficiency at most one.  
 \end{thm}

Theorem~\ref{claw-free} can be restated as $\HH\in 2\text{-}\B\kd$ if and only if $\HH\le\{K_{1,3}\}$.
The main result of this article characterizes the finite family of graphs $\HH\in\B\kd$.

\begin{figure}[tb]
  \centering
\includegraphics[width=6.8cm,height=7.8cm]{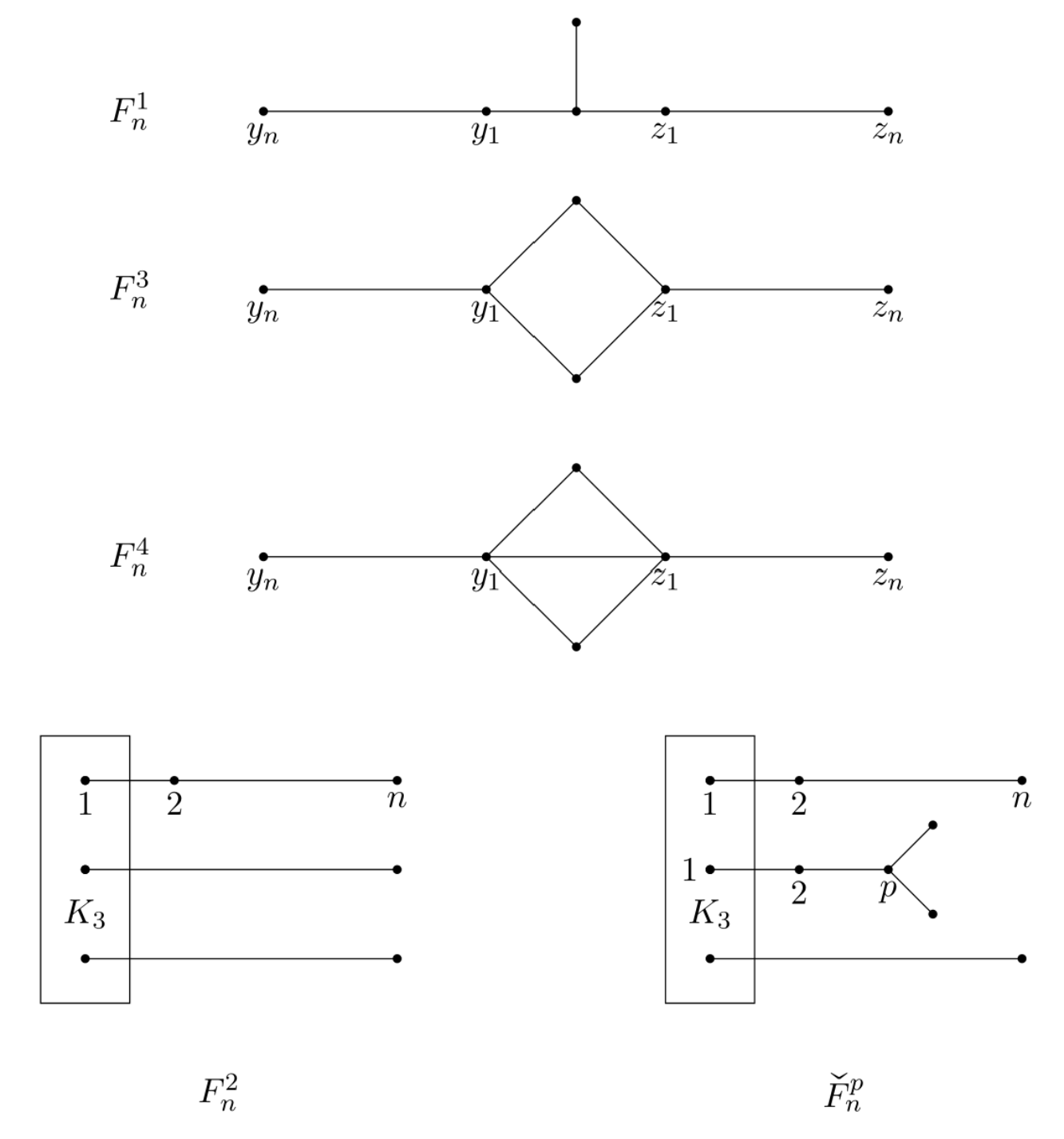}
\caption{Graphs $F_n^1,\,F_n^2,\, F_n^3,\, F_n^4$ and $\fmz F_n^p$.}
  \label{fig 1}
\end{figure}

\begin{figure}[tb]
\centering
\includegraphics[width=6.4cm,height=6.8cm]{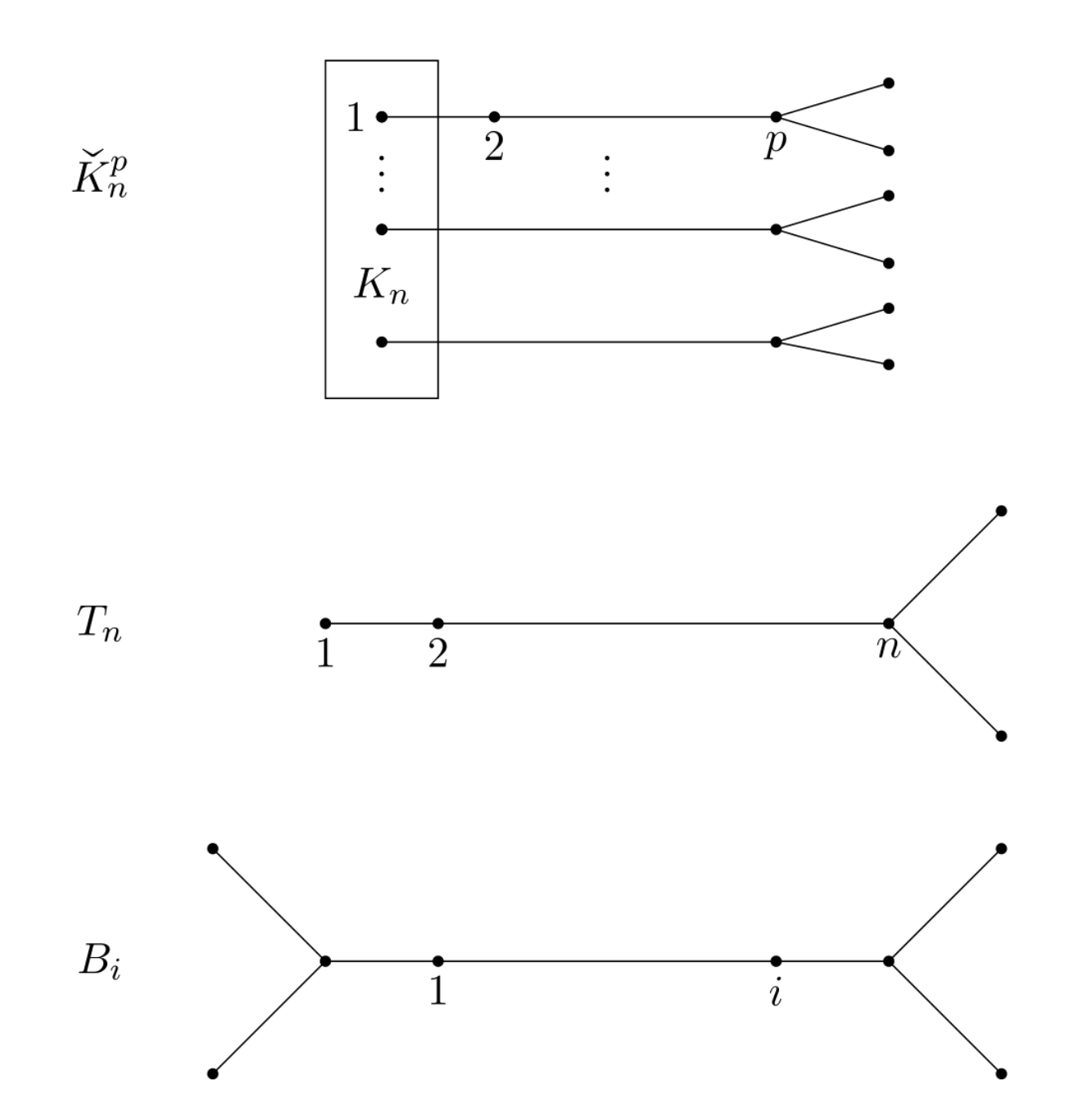}
\caption{Graphs $\fmz K_n^p,\,T_n$ and $B_i$.}
  \label{fig 2}
\end{figure}

In the following constructions, attaching a path $P_n$ at a vertex $v$ of a graph $G$ means that we identify an end of $P_n$ with $v$. We define additional families of graphs that will be used in this article. (See \cref{fig 1,fig 2}.)

\begin{itemize}
   \item $F_n^1$: the graph obtained by attaching  $P_2$ at the center of path $P_{2n+1}$.
  
  \item $F_n^2$: the graph obtained by attaching a copy of $P_n$ at each vertex of triangle $K_3$.

  \item $F_n^3$: the graph obtained by attaching a copy of $P_n$ at each vertex of two fixed nonadjacent vertices of cycle $C_4$.
  
  \item $F_n^4$: the graph obtained by adding an edge to connect the only two vertices of degree $3$ in $F_n^3$. 
 
  \item $T_n$: the graph obtained by attaching two copies of $P_2$ at one end of a path $P_n$. We call this end of $P_n$ {\it the branch vertex} of $T_n$, and the other end of $P_n$ is {\it the end} of $T_n$.
 
\end{itemize}
Note that $T_1$ is actually  $P_3$. The center of $P_3$ is not only the branch vertex but also the end vertex of $T_1$. We continue to define graphs by attaching $T_p$ at a vertex $v$ of a graph $G$, i.e., we identify the end of $T_p$ with $v$.

\begin{itemize}
  \item $B_i \, (i\ge 0)$: the graph obtained by attaching  two copies of $P_2$ at the end of $T_{i+2}$.
  \item $\fmz K_n^p$: the graph obtained by attaching a copy of $T_{p}$ at each vertex of $K_n$.
  \item $\fmz F_n^p$: the graph obtained  by attaching a $T_p$ at one vertex of $K_3$,  and attaching a copy of $P_n$ at each of the remaining vertices of $K_3$. 
\end{itemize}

Now, we state our main theorem.
\begin{thm}\label{main}
Let $\HH$ be a finite family of connected graphs. Then $\HH\in \B\kd$ if and only if 
\[\HH \le \{K_{1,n},\;  T_n, \; \fmz K_n^p : 1\le p\le \frac{n}{2}-1 \}\] or
\[
  \HH \le\{K_{1,n},\; F_n^1,\; F_n^2,\; F_n^3,\; F_n^4,\;\fmz{F}^p_n,\; \fmz {K}_n^p: 1\le p\le n-2\} 
\]     
for some even integer $n\ge 4$.
\end{thm}

For any integer $d\ge 2$, let $d$-$B'$-$\kd$ consist of inclusion-minimal elements of $d$-$B$-$\kd$, i.e. forbidden sets  $\HH\in d$-$B$-$\kd$ while $\HH'\notin d\text{-} \B\kd$ for any proper subset $\HH'\subset \HH$.
According to Theorem~\ref{claw-free}, $\mathcal{H}\in 2$-$B'$-$\kd$ if and only if $\mathcal{H}=\{K_{1,3}\}$.
Plummer and Saito \cite{ps} further showed that, given a graph $H$, for any integer $d\ge 2$,  $\{H\}\in$ $d$-$B'$-$\kd$ if and only if $H\prec K_{1,3}$.
Fujita, Kawarabayashi, Lucchesi, Ota, Plummer, and Saito~\cite{fklops} asked the following question.

\begin{quest}\label{Fujita-Q}
For any integer $d\ge 3$, determine $\HH$ with $|\HH|=2$ and $\HH\in d\text{-}B'\text{-}\kd$.
\end{quest}

Furthermore, they proved the following extension of \cref{claw-free}.
\begin{thm}[\cite{fklops}]\label{bone}
Let $G$ be a connected $\{K_{1,s}, B_i:i\ge 0\}$-free graph where $s\ge 4$. Then $\kd(G)\le s-2$. 
\end{thm}

As an application of \cref{main}, we answer Question~\ref{Fujita-Q}. 

\begin{thm}\label{pair}
	Let $\HH= \{H_1,H_2\}$ and $d\ge 3$ be an integer. 
	Then, $\HH \in d\text{-}B'\text{-}\kd$ if and only if $\HH= \{K_{1,s}, T_3\}$ or $\{K_{1,s},P_4\}$ for some $s$ with $4\le s\le d+1$.
\end{thm}   

\begin{proof}
Note that, for $\HH= \{K_{1,s}, T_3\}$ or $\{K_{1,s},P_4\}$, $\HH\le \{K_{1,s}, B_i:i\ge 0\}$. Thus, by \cref{bone}, we have $\kd(G)\le s-2< d$ since $4\le s\le d+1$, i.e. $\HH\in d\text{-}B\text{-}\kd$. Clearly, $\HH$ is inclusion-minimal. Thus $\HH\in d\text{-}B'\text{-}\kd$.
	
Now suppose $\HH \in d\text{-}B'\text{-}\kd$. Let $\HH_1=\{K_{1,n},\;  T_n, \; \fmz K_n^p : 1\le p\le \frac{n}{2}-1 \}$ and $\HH_2=\{K_{1,n},\; F_n^1,\; F_n^2,\; F_n^3,\; F_n^4,\; \fmz {K}_n^p, \;\fmz{F}^p_n: 1\le p\le n-2\}$. Thus, $\HH_1, \HH_2 \le  \{K_{1,n}, \, F_n^1,\, \fmz K_n^1 \}$.
Since $\HH \in d\text{-}B'\text{-}\kd \subset \B\kd $, we have $\HH\le \{K_{1,n},\, F_n^1,\, \fmz K_n^1 \}$. Without loss of generality, assume $H_1\prec K_{1,n}$. By \cref{claw-free} and the minimality of $\HH$,  we obtain $H_1=K_{1,s}$ for some $s\ge 4$. Thus, $H_1$ is not an induced subgraph of  $F_n^1$ or $\fmz K_n^1$. Hence $\{H_2\}\le \{F_n^1, \fmz K_n^1\}$. However, $H_2\prec F_n^1$ implies that $H_2$ is triangle-free.
Hence,  $H_2$ contains at most two vertices of the largest clique of $\fmz K_n^1$. Therefore, $H_2\prec B_0$. Since $B_0\nprec F_n^1$, we have $H_2\prec T_3$. By the minimality of $\HH$, we have $H_2=$ $T_3$ or $P_4$. Therefore, we have $\HH=\{K_{1,s}, T_3\}$ or $\{K_{1,s}, P_4\}$. 
$s\le d+1$ as $\kd(K_{1,s-1})=s-2<d$ and $K_{1,s-1}$ is $\HH$-free. 
\end{proof}

The rest of this article is arranged as follows. We present some preliminaries in Section 2 and prove Theorem~\ref{main} in Section 3. We provide some discussion and remarks in the last section.

\section{Preliminaries}

We first construct several graphs with large deficiencies. 
Let $s,t$ be positive integers,  and let $Q_i=u_i^1u_i^2\dots u_i^t$, $i\in [s+1]$ be $s+1$ pairwise vertex disjoint paths.
Define
\begin{itemize}
 \item  $H_{s,t}^1$: the graph obtained from the union of paths $Q_i, i\in [s+1]$ by adding $2s$ vertices $\{v_i, w_i: i\in [s]\}$ and $3s$ edges 
$ \{v_iw_i, \;v_iu_i^t,\; v_iu_{i+1}^1: i\in[s]\}$.

 \item $H_{s,t}^3$: the graph obtained from the union of paths $Q_i, i\in [s+1]$ by adding $2s+2$ vertices $\{x,y\}\cup\{v_i,w_i : i\in [s]\}$ and $4s+2$ edges 
 $\{xu_1^1, u_{s+1}^ty\}\cup\{v_iu_i^t, v_iu_{i+1}^1, w_iu_i^t, w_iu_{i+1}^1 : i\in [s]\}$.

 \item $H_{s,t}^4$: the graph obtained from $H_{s,t}^3$ by adding $s$ edges $\{u_i^tu_{i+1}^1: i\in [s]\}$.

 \item $\fmz H_{s,t}^p$: the graph obtained from the union of paths $Q_i, i\in [s+1]$ and $R_i\cong T_p, i\in [s]$ with end $v_i$ by adding $3s$ edges $\{u_i^tu_{i+1}^1, u_i^tv_i, u_{i+1}^1v_i : i\in [s]\}$.
\end{itemize}
For convience, let $H_{0,t}^1=P_t$ and $H_{0,t}^3=H_{0,t}^4=P_{t+2}$. 

\begin{lem}\label{H}
(1) \begin{equation*}
\kd (\fmz H_{s,t}^p)=
\begin{cases}
 s\qquad & \mbox{if } t(s+1)+s(p+2) \mbox{ is even,}\\
                      s+1 & \mbox{otherwise.}
\end{cases} 
\end{equation*}

(2) \begin{equation*}
	\kd (\fmz K_{n}^p)=
	\begin{cases}
		n\qquad &\mbox{if } n(p+1) \mbox{ is even,}\\
		n+1 & \mbox{otherwise.}
	\end{cases} 
\end{equation*}
\end{lem}

\begin{proof}
Recall the definition of $\fmz H_{s,t}^p$ (resp. $\fmz K_{n}^p$), for any $i\in [s]$ (resp. $i\in[n]$), each $R_i(\cong T_p)$  contains a branch vertex that is adjacent to two  leaves. Given a maximum matching $M$ of $\fmz H_{s,t}^p$ (or $\fmz K_{n}^p$), $M$ misses at least one of these two leaves. Thus $M$ misses at least $s$ (resp. $n$) vertices of $\fmz H_{s,t}^p$ (resp. $\fmz K_{n}^p$). Therefore, $\kd (\fmz H_{s,t}^p)\ge s$ (resp. $\kd(\fmz K_{n}^p)\ge n$).

Let $H$ be the graph obtained by deleting a leaf adjacent to every branch vertex from $\fmz H_{s,t}^p$ (resp. $\fmz K_{n}^p$). Then $H$ is connected and claw-free. By \cref{claw-free}, $\kd(H)\le 1$. Hence $\kd (\fmz H_{s,t}^p) \le s+1$ (resp. $\kd(\fmz K_{n}^p)\le n+1$).

Since $\kd(G)\equiv |V(G)|\pmod 2$ for any graph $G$, the proof is complete. 
\end{proof}

Next we handle $H_{s,t}^j: j\in \{1,3,4\}$. 
When $t$ is even, it can be easily checked that $\kd(H_{s,t}^1)= \kd(H_{s,t}^3)=\kd(H_{s,t}^4)=0$, which is not what we want. When $t$ is odd, we show that $H_{s,t}^j$ has large deficiency for $j\in\{1,3,4\}$.
\begin{lem}\label{H^1}
If $t$ is odd, then $\kd(H_{s,t}^1)=s+1$. 
\end{lem}

\begin{proof}
Let $H=H_{s,t}^1-\{w_1, w_2, \cdots, w_s\}$. Then $H$ is a path, thus $\kd(H)\le 1$. As a maximum matching of $H$  is also a matching of $H_{s,t}^1$, we have  $\kd(H_{s,t}^1)\le s+1$.

Now we prove that  $\kd(H_{s,t}^1) \ge s+1$ by induction on $s$. For the base case $s=0$, $H_{0,t}^1=P_{t}$ has odd vertices. Thus $\kd(H_{0,t}^1)\ge 1$.
Assume that $s\ge 1$.  Let $M$ be a maximum matching  of $H_{s,t}^1$. If $v_1\not \in V(M)$, then $w_1$ must not be in $V(M)$. Thus $M\cup \{v_1w_1\}$ is a larger matching, a contradiction. 
Therefore, $v_1\in V(M)$.
If $\{v_1u_1^t\}\in M$ (or $\{v_1u_2^1\}\in M$), then $w_1$ is missed by $M$. Hence $M'=M\cup \{v_1w_1\} \setminus \{v_1u_1^t\}$ (or $M'=M\cup \{v_1w_1\} \setminus \{v_1u_2^1\}$) is still a maximum matching.
Thus we may assume that $\{v_1w_1\}\in M$. Then $H'=H_{s,t}^1-\{v_1,w_1\}$ consists of two components, $H_{0,t}^1$ and $H_{s-1,t}^1$. By induction hypothesis, $\kd(H_{0,t}^1)\ge 1$ and $\kd(H_{s-1,t}^1)\ge s$. Hence $M$ misses at least $1$ vertex in $H_{0,t}^1$ and $s$ vertices in $H_{s-1,t}^1$. Therefore, $\kd(H_{s,t}^1)\ge s+1$.
\end{proof}

\begin{lem}\label{H^3}
If $t$ is odd, then $\kd(H_{s,t}^3)=\kd(H_{s,t}^4)=s+1$. 
\end{lem}

\begin{proof}
Since $H_{s,t}^3$ is a spanning  subgraph of $H_{s,t}^4$, we have $\kd(H_{s,t}^4)\le \kd(H_{s,t}^3)$. 
Let $H=H_{s,t}^3-\{w_1, w_2, \cdots, w_s\}$. Then $H$ is a path too. Thus $\kd(H)\le 1$. As a maximum matching of $H$  is also a matching of $H_{s,t}^3$, we have  $\kd(H_{s,t}^3)\le s+1$ too.

Now we prove that  $\kd(H_{s,t}^4) \ge s+1$ by induction on $s$. For the base case $s=0$, $H_{0,t}^4=P_{t+2}$. Thus $\kd(H_{0,t}^4)\ge 1$ as $t$ is odd.
Assume that $s\ge 1$.  Let $M$ be a maximum matching  of $H_{s,t}^4$.  
If $u_1^t\notin V(M)$, then at least one of $v_1$ and $w_1$ cannot be saturated by $M$ since both can be matched with only $u_2^1$. By symmetry,  we may assume $v_1\notin V(M)$. Then $M\cup\{u_1^tv_1\}$ is a larger matching of $H_{s,t}^4$, a contradiction.
Now assume $u_1^t\in V(M)$. If $u_1^tu_2^1\in M$, then both $v_1$ and $w_1$ are missed by $M$. Therefore,
$(M\setminus\{u_1^tu_2^1\})\cup \{u_1^tv_1,\; u_2^1w_1\}$ is a larger matching, a contradiction too. 
If $u_{1}^tu_{1}^{t-1}\in M$ and $v_1$ is not saturated by $M$, then $M'=M\setminus \{u_1^tu_1^{t-1}\}\cup \{u_1^tv_1\}$ is still a maximum matching.

Thus we may always assume $u_{1}^tv_1\in M$ by symmetry of $v_1$ and $w_1$ in $H_{s,t}^4$. 
Let $H'=H_{s,t}^4-\{u_1^t,v_1\}$, then $H'$ consists of two components, $H_{0, t-2}^4$ and $H_{s-1,t}^4$. By induction hypothesis, $\kd (H_{0, t-2}^4)\ge 1$ and $\kd(H_{s-1,t}^4)\ge s$.  Hence $M$ misses at least one vertex in $H_{0,t}^4$ and at least $s$ vertices in $H_{s-1,t}^4$. Therefore, $\kd(H_{s,t}^4)\ge s+1$. The proof is complete.
\end{proof}

Let $X$ be a set of vertices in a graph $G$. Denote by $\delta_G(X)=\min\{d_G(x) : x\in X\}$ and $\Delta_G(X)=\max\{d_G(x) : x\in X\}$. For two vertex sets $X$ and  $Y$ in $G$, we say that $X$ {\em dominates} $Y$ if $Y\subseteq N_G(X)$. 
Recall that an induced matching $M$ of $G$ is an induced 1-regular subgraph of $G$. 
The {\em induced matching number} of $G$ is defined as $$\nu' (G)=\max\{|M| : \text{$M$ is an induced matching of $G$} \}.$$ 
Next we present a technical lemma that has also been used implicitly in \cite{cf22,f}.
\begin{lem}\label{induced matching}
	Let $G$ be a bipartite graph with vertex partition sets $X$ and $Y$. Suppose $\delta(X)\ge 1$ and $\Delta(Y)\le n$.
	If  $|X|\ge n(p-1)+1 $, then $\nu' (G)\ge p$. Moreover, the lower bound of $|X|$ is tight.
\end{lem}
\begin{proof}
Since $\delta(X)\ge 1$, $Y$ dominates $X$. Take a minimal subset $Y'\subset Y$ dominating $X$, i.e. $N_G(Y')=X$ and for any proper subset $Y''$ of $Y'$, $N_G(Y'')\not =X$. Therefore, $Y'$ is irredundant, i.e., every $y\in Y'$ has a {\em private} neighbor $x_y\in X$ (i.e., $N_G(x_y)\cap Y'=\{y\}$). Since $\Delta(Y)\le n$, every $y\in Y'$ dominates at most $n$ vertices in $X$. Since $|X|\ge n(p-1)+1$, we have  $|Y'|\ge p$. Hence, we can choose $p$ vertices in $Y'$ and a private neighbor for every one; they together form an induced matching $M_p$. Therefore, $\nu' (G)\ge p$.  

Let $G=(X\cup Y,E)$ be a bipartite graph with  $X=\{x_{ij}:i\in [p-1],\, j\in [n]\}$, $Y=\{y_i:i \in [p-1]\}$, and $E=\{(x_{ij},y_i):i \in [p-1]\}$. Then we have  $\delta(X)=1$, $\Delta(Y)=n$, but $\nu' (G)= p-1$. 
\end{proof}

For two disjoint vertex sets $X$ and $Y$ of a graph $G$, write $G(X, Y)$ for the bipartite subgraph induced by the edges between $X$ and $Y$ in $G$. The following result  is  important in this article and has its own meaning on this subject. 

\begin{thm}\label{finite level}
Let $n$ and $N$ be two integers greater than $3$. Let $G$ be a connected  $\{ K_{1,n}, \fmz{K}_n^p: p\ge 1\}$-free graph. If $\diam(G)\le N$, then there is an integer $f(n, N)$ such that $\kd(G)\le f(n, N)$.  
\end{thm}

\begin{proof}
Fix a vertex $x_0\in V(G)$, for $i\ge 0$, define 
$$
X_i=\{v\in V(G):\dist(x_0,v)=i\}.
$$ 
Since $\diam(G)\le N$, we may assume that $X_N$ is the last level. Then $V(G)=\bigcup_{i=0}^N X_i$. 
Now we construct a matching level by level. First, find a maximum matching $M_N$ in $G[X_N]$ and let $Y_N=X_N\setminus V(M_N)$. Then $Y_N$ is a stable set. Next, find a maximum matching $M'_N$ in $G(X_{N-1},Y_N)$ such that $G_{N-1}=G-V(M_N)-V(M'_N)$ is still connected. Let $Z_N=Y_N-V(M_N')$ and $B_{N-1}=N_{G_{N-1}}(Z_N)$. Let $H_N=G(Z_N, B_{N-1})$.  
	
\begin{claim}\label{snail horn}
For every vertex $u\in B_{N-1}$, there exist $x$ and $y$ in $Z_N$ such that $N_{G_{N-1}}(x)=N_{G_{N-1}}(y)=\{u\}$, i.e., $u$ has at least two private neighbors in $Z_N$.
\end{claim}  
\begin{proof}
Since $u\in B_{N-1}=N_{G_{N-1}}(Z_N)$, there is a vertex $z\in Z_N$ adjacent to $u$. Then $G_{N-1}-\{u,z\}$ is disconnected; otherwise, $M'_N\cup \{uz\}$ will be a larger matching in $G(X_{N-1},Y_N)$ such that $G-V(M_N)-V(M'\cup \{uz\})$ is still connected. Hence there is a vertex $x\in Z_N$ such that $N_{G_{N-1}}(x)=\{u\}$. For  the same reason, $G_{N-1}-\{x,u\}$ is disconnected. Thus we can find another vertex $y\in Z_N$ (possibly $z$) such that $N_{G_{N-1}}(y)=\{u\}$.   
\end{proof}

	\begin{claim}\label{DeltaH}
	If $N\ge 2$, then	$\Delta_{H_N}(B_{N-1})\le n-2$.	
	\end{claim}
	\begin{proof}
		Suppose to the contrary that $\Delta_{H_N}(B_{N-1})\ge n-1$. Then the neighbors of a vertex $u\in B_{N-1}$ with $d_{H_N}(u)=n-1$ together with a neighbor of $u$ in $X_{N-2}$ induces a star $K_{1,n}$, a contradiction.
	\end{proof}

	Set $\alpha_1=n$ and $\beta_1=R(n,n)$. For $i\ge 1$, recursively define $$
	\alpha_{i+1}=(n-2)(\beta_i-1)+1,\quad \beta_{i+1}=R(n,\alpha_{i+1}).
	$$
	
  \begin{claim}
		$|Z_N|\le (\beta_{N-1}-1)(n-2)$.  
	\end{claim} 
	\begin{proof}
Suppose to the contrary that $|Z_N|> (\beta_{N-1}-1)(n-2)$. Then By Claim~\ref{DeltaH}, we have $|B_{N-1}|\ge |Z_N|/(n-2) >\beta_{N-1}-1$, i.e., $|B_{N-1}| \ge\beta_{N-1}$.
By \cref{snail horn}, for each vertex $u\in B_{N-1}$, $u$ has at least two private neighbors $x$ and $y$ in $Z_N$. Recall that $|B_{N-1}|\ge \beta_{N-1}=R(n, \alpha_{N-1})$. By the Ramsey Theorem, $G[B_{N-1}]$ contains a clique $K_n$ or a stable set of order $\alpha_{N-1}$. If $G[B_{N-1}]$ contains a clique $K_n$, then $G[V(K_n)\cup N_{H_{N}}(V(K_n))]$ must contain an induced $\fmz K_n^1$, a contradiction.

Now, assume that there is a stable set $A_{N-1}$ in $G[B_{N-1}]$ with $|A_{N-1}|\ge \alpha_{N-1}$. Set $H_{N-1}=G(A_{N-1}, X_{N-2})$.
		Clearly, $\delta_{H_{N-1}}(A_{N-1})\ge 1$. For the same reason as in Claim~\ref{DeltaH}, we have $\Delta_{H_{N-1}}(X_{N-2})\le n-2$.
		By \cref{induced matching},  $H_{N-1}$ has an induced matching $M_{N-1}$ with $|M_{N-1}|\ge \beta_{N-2}$. 
		Let $B_{N-2}=V(M_{N-1})\cap X_{N-2}$. Then $|B_{N-2}|=|M_{N-1}|\ge \beta_{N-2}=R(n,\alpha_{N-2})$.
		Note that each vertex of $V(M_{N-1})\cap A_{N-1}$ has at least two private neighbors in $Z_N$. Hence, if $G[B_{N-2}]$ contains a clique $K_n$, then $G$ must contain an induced $\fmz K_n^2$, a contradiction.  Thus there is a stable set $A_{N-2}$ in $G[B_{N-2}]$ with $|A_{N-2}|\ge \alpha_{N-2}$. 
		Set $H_{N-2}=G(A_{N-2}, X_{N-3})$. With the same discussion as for $H_{N-1}$, $H_{N-2}$ has  an induced matching $M_{N-2}$ with $|M_{N-2}|=|B_{N-3}|\ge \beta_{N-3}$ and a stable set $A_{N-3}\subseteq B_{N-3}$ with $|A_{N-3}|\ge \alpha_{N-3}$, where $B_{N-3}=V(M_{N-2})\cap X_{N-3}$. 
	The above process is continued until we find a stable set $A_1\subseteq B_1\subseteq X_1$ with $|A_1|\ge\alpha_1=n$. However, $\{x_0\}\cup A_1$ induces a copy of $K_{1,n}$, a contradiction. This claim holds.
	
	\end{proof}

	Recall that every vertex of $B_{N-1}$ has a private neighbor in $Z_N$. Hence $H_N$ has an induced matching $M_N''$ saturating all vertices of $B_{N-1}$. Let $Z'_N=Z_N\setminus V(M_N'')$. Then 
	$$
	|Z'_N|=|Z_N|-|B_{N-1}|\le \left(1-\frac{1}{n-2}\right)|Z_N|\le (n-3)(\beta_{N-1}-1).
	$$ 
	Therefore, the matching $M_N\cup M_N'\cup M_N''$ saturates the vertices in $X_N$  except for $|Z_N'|\le (n-3)(\beta_{N-1}-1)$ vertices.
	Let $G'_{N-1}$ be the graph obtained from $G$ by deleting all vertices in $V(M_N)\cup V(M_N')\cup V(M_N'')\cup Z'_N$ and update $X_{N-1}$ by $X_{N-1}-(V(M_N')\cup V(M_N''))$.
	Then $V(G'_{N-1})=\bigcup_{i=0}^{N-1}X_i$ and $G'_{N-1}$ is still connected.
	
	We can continue the above process to $G'_{i}$ for each $i=N-1, N-2, \dots, 2$ such that we can construct the matching $M_i\cup M_i'\cup M_i''$ in $G'_{i}$ that  saturates vertices in $X_i$ except for $|Z_i'|\le (n-3)(\beta_{i-1}-1)$ vertices. Finally, let $G'_1=G'_2-\left(V(M_2)\cup V(M_2')\cup V(M_2'')\cup Z_2'\right)$ and update $X_1$ by $X_1-(V(M_2')\cup V(M_2''))$. 

	Find a maximum matching $M_1$ in $G[X_1]$ and let $Y_1=X_1-V(M_1)$. Thus $Y_1$ is a stable set dominated by $x_0$. Since $G$ is $K_{1,n}$-free, we have $|Y_1|\le n-1$. Let $M_1'$ be a maximum matching of $G[\{x_0\}\cup Y_1]$. Set $M_1''=\emptyset$.

Finally, we have a matching $\bigcup_{i=1}^{N} (M_i\cup M_i'\cup M_i'')$ of $G$ missing at most $$
	n-2+(n-3)\sum_{i=1}^{N-1}(\beta_i-1)=f(n,N)
	$$
	vertices. Thus $\kd(G)\le f(n,N)$. 
\end{proof}

\section{Proof of Theorem~\ref{main}}
\begin{proof}[\bf The ``only if" part:]  Suppose that $\HH\in \B\kd$, i.e., there is a constant $c=c(\HH)$ such that every connected $\HH$-free graph $G$ satisfies $\kd(G)< c$. Since $\HH$ is a finite family, the value $h=\max \{|V(H)|:H\in \HH\}$ is well-defined. 

Choose $s>c$, and $p>h$ as integers and $t>h$ as an odd integer. By Lemma~\ref{H}, we have $\kd(\fmz H_{s,t}^p)\ge s>c$. Hence $\fmz H_{s,t}^p$ is not $\HH$-free, i.e., there is an $H\in \HH$ such that   $H\prec\fmz H_{s,t}^p$.  Since  $t>h\ge|V(H)|$ and $p>h\ge|V(H)|$,  $H$ contains at most one of the end vertex and the branch vertex of each $R_i: i\in [s]$. Thus there is some even integer $n>h\ge |V(H)|$ (we choose even $n$ just for convenience) such that $H\prec T_n$ or $H\prec F_n^2$. We may assume $n>c+1$ further.

If $H\prec T_n$, then $\HH \le \{T_n\}$. Since 
\[
\kd(K_{1,n})=n-1>c,\text{ and } \kd (\fmz K_n^p)\ge n>c,\; \text{for any} \;  p\in\left [\frac{n}{2}-1\right],
\]
$K_{1,n}$ and  $\fmz K_n^p: 1\le p\le \frac{n}{2}-1$ are not $\HH$-free either. Therefore, $\HH \le \{K_{1,n}, T_n,  \fmz K_n^p : 1\le p\le \frac{n}{2}-1 \}$.

For the other case $H\prec F_n^2$, we have $\HH \le \{F_n^2\}$. By \cref{H^1} and \ref{H^3},  $\kd(H_{s,t}^j)= s+1>c$ for $j\in\{1,3,4\}$. Thus $H_{s,t}^j$ is not $\HH$-free for $j \in \{1,3,4\}$. Hence there exists  $H_j\in \HH$ such that  $H_j\prec H_{s,t}^j$. 
Since $|V(H_j)|\le h$, we have  $|\{i\in [s]: V(H_j)\cap \{v_i, w_i\}\neq \emptyset \}|\le 1$. This implies that $H_j\le F_n^j$. 
By Lemma~\ref{H}, $\kd (\fmz H_{s,t}^p)>s>c$ for $1\le p\le n-2$. Thus $\fmz H_{s,t}^p$ contains an induced subgraph $H_p\in \HH$ too. Since $|H_p|\le h< t$, 
$|\{i\in [s]: v_i\in V(H_p) \}|\le 1$ too.
This implies that $H_p\prec \fmz F_n^p$. 
Note that both $K_{1,n}$ and $\fmz K_n^p: 1\le p\le n-2$ are not $\HH$-free. In summary, we have 
\[\HH \le\{K_{1,n}, F_n^1, F_n^2, F_n^3, F_n^4, \fmz {K}_n^p, \fmz{F}^p_n: 1\le p\le n-2\} .\] 
The proof of the ``only if'' part is completed.
\end{proof}

\begin{proof}[\bf The ``if'' part:]
Either  
$$\HH \le \{K_{1,n},\; T_n,\;  \fmz K_n^p : 1\le p\le \frac{n}{2}-1 \}$$
or  
$$\HH \le\{K_{1,n},\; F_n^1,\; F_n^2,\; F_n^3,\; F_n^4,\; \fmz {F}_n^p,\; \fmz{K}^p_n: 1\le p\le n-2\}$$
implies that every $\HH$-free graph $G$ must be
$\{K_{1,n},\; F_n^1,\; F_n^3,\; F_n^4,\;\fmz F_n^p,\;\fmz K_n^p: p\ge 1\}$-free.

Take a longest induced path $P=u_1u_2\dots u_m$ of $G$. We may assume $m\ge n^2-n-1$, otherwise, by \cref{finite level} the proof is completed since $\diam(G)\le m-1<n^2$. 
We consider the following decomposition of $G$. (See \cref{fig 3}.) 
\begin{figure}[tb]
  \centering
\includegraphics[width=8cm, height=6cm]{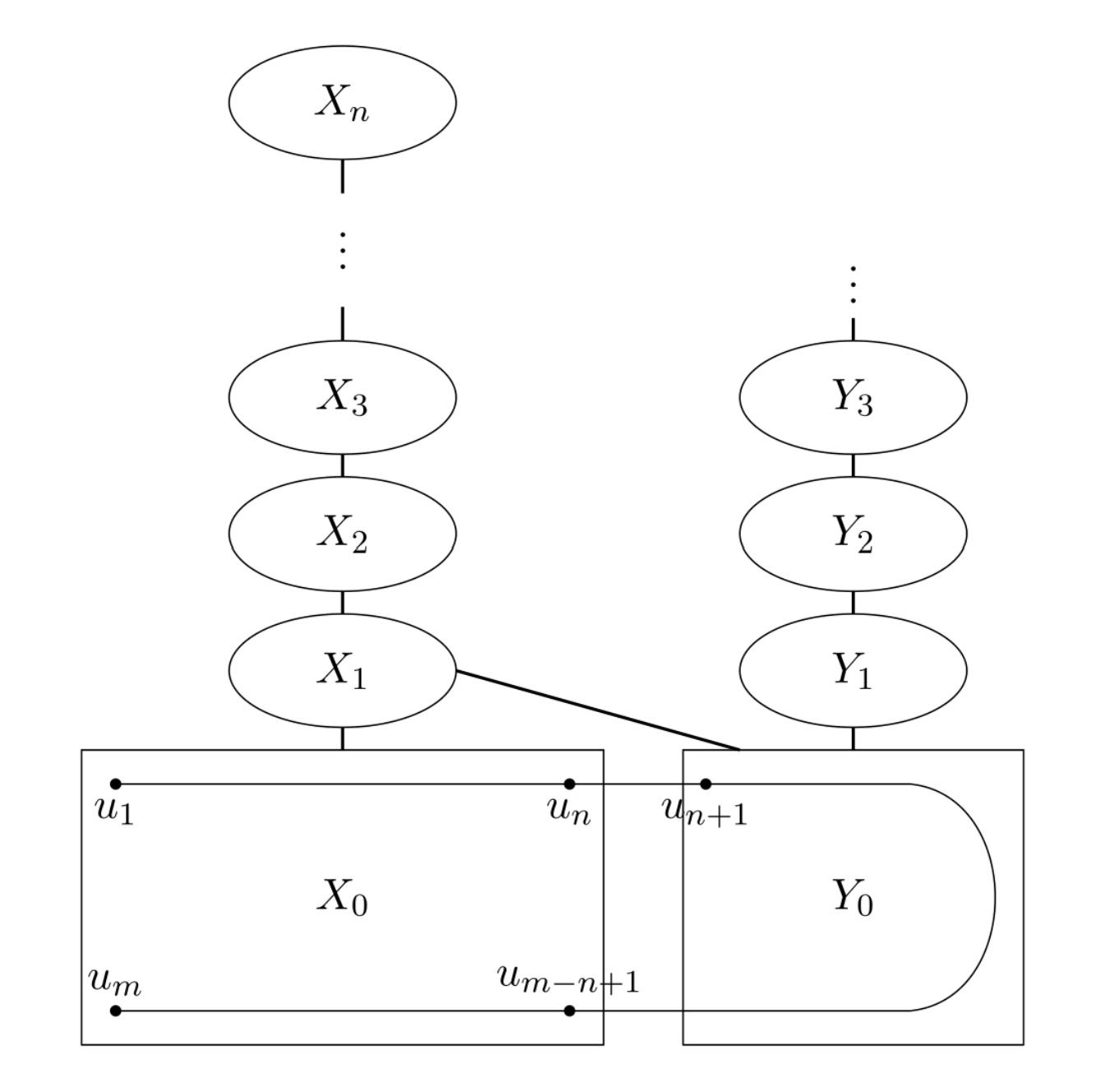}

  \caption{The decomposition of $G$.}
  \label{fig 3}
\end{figure}

Let \[X_0=\{u_i: 1\le i\le n \text{ or } m-n+1\le i\le m\}.\]
Let $Y_0=\{u_i: n+1\le i\le m-n\}.$ 
Set

 \[X_1=N_G(X_0)-V(P),\text{ and } Y_1=N_G(Y_0)-(X_0\cup X_1).\] 
For $i\ge 2$, recursively define 
\[X_i=N_G(X_{i-1})-(Y_0\cup Y_1)-\bigcup_{j=1}^{i-2}X_j.\] 

\begin{claim}\label{CL: X_2n0+1}
 $X_{n+1}=\emptyset$. 
\end{claim}
\proof
Suppose to the contrary that there exists a vertex $x_{n+1}$ in $X_{n+1}$.
Then we can recursively trace back to its ancestors $x_{n+1-i}\in N_G(x_{n+2-i})\cap X_{n+1-i}$ for each $i\in [n+1]$. 
Assume $x_{0}=u_{k}$ for some $k$ with $1\le k\le n$ or $ m-n+1\le k\le m$.
By symmetry, we may assume $k\in [n]$ and that such a $k$ is as large as possible in $[n]$.
Consider the vertices in $N_G(x_{1})\cap V(P)$. 
Let $j_{1}=\min\{j\in [m]:x_{1}u_{j}\in E(G)\}$.
For $p\geq 2$, we recursively define  \[j_{p}=\min\{j:j_{p-1}+2\leq j\leq m,~x_{1}u_{j}\in E(G)\},\]
and stop when \[\{j:j_{p}+2\leq j\leq m,~x_{1}u_{j}\in E(G)\}= \emptyset.\]
Let $S=\{u_{j_{p}}:p\geq 1\}$, and set $s=|S|$.
Since $j_{p}\geq j_{p-1}+2$ and $P=u_1u_2\dots u_m$ is an induced path, we have $S$ is a stable set of $G$.
Thus $\{x_{1},x_{2}\}\cup S$ induces a copy of $K_{1,s+1}$ with center $x_1$ in $G$.
Since $G$ is $K_{1,n}$-free, we have $s+1\leq n-1$.

We claim that $s\ge 2$. If $s=1$, then $\alpha(N(x_1)\cap V(P))=1$. Since $x_1u_k\in E(G)$, $x_1u_i\notin E(G)$ for any $i\ge k+2$. If $x_1u_{k+1}\in E(G)$, then $x_{n+1}x_n\dots x_1u_{k+1}u_{k+2}\dots u_m$ is an induced path with $n+m+1-k>m$ vertices, a contradiction to the longest $P$. On the other hand, if $x_1u_{k+1}\notin E(G)$, then $x_{n+1}x_n\dots x_1u_ku_{k+1}u_{k+2}\dots u_m$ is an induced path with $n+m+2-k>m$ vertices, which is also a contradiction.

Let $j^{*}=\max\{j\in [m]:\, x_{1}u_{j}\in E(G)\}$.
Note that $j^{*}\in \{j_{s},j_{s}+1\}$. 
Let $Q_{1}=u_{1}u_{2}\cdots u_{j_{1}}$,  $Q_{p}=u_{j_{p-1}+2}u_{j_{p-1}+3}\cdots u_{j_{p}}$ for $p\in [2,s]$, and $Q_{s+1}=u_{j^{*}}u_{j^{*}+1}\cdots u_{m}$.
Then \[V(P)- \bigcup _{p=1}^{s+1} V(Q_{p})=\{u_{j_{p}+1}: p\in [s-1]\},\] and hence
\begin{align*}
m= |V(P)|&= \left| V(P)- \bigcup_{p=1}^{s+1} V(Q_{p}) \right|+ \left|\bigcup_{p=1}^{s+1} V(Q_p) \right|\\
&\leq (s-1)+\sum_{p=1}^{s+1} |V(Q_{p})|\\
&= \sum_{p=1}^{s+1} (|V(Q_{p})|+1)-2.
\end{align*}
This implies that 
\[
\sum_{p=1}^{s+1}(|V(Q_{p})|+1)\ge m+2 \geq n^{2}-n+1.
\]
If $|V(Q_{p})|\leq n-1$ for all $p\in [s+1]$, then \[
n^{2}-n+1\leq \sum_{p=1}^{s+1}(|V(Q_{p})|+1) \leq (s+1)n\leq (n-1)n,
\] 
a contradiction. Thus there is some $q\in [s+1]$ such that  $|V(Q_{q})|\geq n$.
Note that $|N_G(x_{1})\cap V(Q_{q})|=1$.
Write $N_G(x_{1})\cap V(Q_{q})=\{u\}$, i.e., $u=u_{j_q}$ if $q\neq s+1$, and $u=u_{j^{*}}\in \{u_{j_{s}}, u_{j_s+1}\}$ if $q=s+1$.
Since $|S|\geq 2$, we can choose anther vertex $v\in S$ such that $N_G(v)\cap V(Q_{q})=\emptyset$ (e.g., choose $v\in S- \{u\}$  if $q\neq s+1$, and, otherwise, 
let $v=u_{j_{1}}$).
Note that $u$ is an end vertex of $Q_{q}$, $|V(Q_{q})|\ge n$. 
We have a copy of $F^1_n\prec G[\{v,x_{1},x_{2},\ldots ,x_{n+1}\}\cup V(Q_q)]$, a contradiction.
\qed

Now set \[
X=\bigcup_{i=0}^{n}X_i.
\]
Recall that $Y_0=\{u_i: n+1\le i\le  m-n\}$, and $Y_1=N_G(Y_0)-X$. For $i\ge 2$, recursively define \[
Y_i=N_G(Y_{i-1})-X-\bigcup_{j=0}^{i-2}Y_j.
\] 
Suppose $N$ is the largest index with $Y_N\neq \emptyset$ and $Y_{N+1}=\emptyset$.
Set \[Y=\bigcup_{i=0}^{N} Y_i.\]
Then $V(G)=X\cup Y$.
Take vertices $v,w\in Y_1$ (if any) such that 
$$i=\min \{j: vu_j\in E(G) \}= \min \{j: wu_j\in E(G) \}.$$
Let $$j=\max \{k: vu_k\in E(G) \}, \text{ and } j'=\max \{k: wu_k\in E(G) \}.$$ 
By symmetry, we may assume $j\ge j'$.
Note that $Y_1\cap N_G(X_0)=\emptyset$. Hence $ n+1\le i\le j'\le j\le m-n$. 

\begin{claim}\label{i+1}
$vu_{i+1}\in E(G)$.
\end{claim}
\begin{proof}
Suppose to the contrary that $vu_{i+1}\notin E(G)$.  Thus $j\neq i+1$. 
If $j=i$, then

  \[\{u_{i-n}, u_{i-n+1}, \dots, u_{i-1}, u_i, v, u_{i+1}, u_{i+2}, \dots, u_{i+n}\}\]
induces a copy of $F^1_n$ in $G$, a contradiction.
If $j=i+2$, then \[
  \{u_{i-n+1}, u_{i-n+2},\dots, u_{i}, v, u_{i+1},u_{i+2},u_{i+3},\dots ,u_{i+n+1}\}
\]
induces a copy of $F^3_n$ in $G$, a contradiction again.
Now assume $j>i+2$. Then 
$$
\{u_{i-n}, u_{i-n+1}, \dots, u_{i-1}, u_i, u_{i+1}, v, u_{j}, u_{j+1},\ldots ,u_{j+n-2}\}
$$
induces a copy of $F_n^1$, a contradiction. This completes the proof of the claim. 
\end{proof}

\begin{claim}\label{child}
  If there exists $r\in N(v)\cap Y_2$, then $N_G(v)\cap Y_0=\{u_i,u_{i+1} \}$.
\end{claim}
\begin{proof}
 By \cref{i+1}, $j\ge i+1$. It is sufficient to show that $j=i+1$. Suppose to the contrary that $u_i$ and $u_j$ are not adjacent in $P$. Then  
 \[
    \{u_{i-n+1}, u_{i-n}, \dots,  u_{i-1},  u_i, v,r, u_j, u_{j+1}, \ldots ,u_{j+n-1} \}
  \] 
induces a copy of $F_n^1$ in $G$, a contradiction.
\end{proof}

\begin{claim}\label{pq}
$vw\in E(G)$.  
\end{claim}
\proof 
Suppose to the contrary that $vw\not \in E(G)$.
If $j=j'$, then
$$
\{u_{i-n+1}, u_{i-n+2}, \dots, u_{i}, v,w,u_{j},u_{j+1},\ldots ,u_{j+n-1}\}
$$
induces a copy of $F_n^3$ or $F_n^4$ in $G$ according to whether $j=i+1$ or not, which is a contradiction.
Thus $j>j'$. Therefore,  $wu_{j}\notin E(G)$. Then
$$
\{ u_{i-n},u_{i-n+1},\ldots, u_{i-1}, u_i,w, v,u_{j},u_{j+1},\ldots ,u_{j+n-2}\}
$$
induces a copy of $F_n^1$ in $G$, a contradiction too.\qed\bigskip

Since $G$ is $\fmz F_n^p$-free for any $p\ge 1$, 
we know that for any $y\in Y_i$, $i\ge 1$, the children of $y$ in $Y_{i+1}$ are pairwise adjacent. Otherwise, we can construct a  copy of  $\fmz F_n^p$ for some $p\ge 1$ as an induced subgraph in $G$ according to \cref{child}. In other words, $N_G(y)\cap Y_{i+1}$ induces a clique in $G$ for any $y\in Y_i$, $i\ge 1$.

Now we construct a matching step by step from the last level $Y_N$.
First, we find a maximum matching $M_N$ in $G[Y_N]$. Then $Z_N:=Y_N\setminus V(M_N)$ is a stable set. Hence for every vertex $z\in Z_N$, there exists a $p_z\in Y_{N-1}$ such that $N(p_z)\cap Z_N=\{z\}$ as long as $N\ge 2$. Therefore, $M_N'=\{zp_z : z\in Z_N \text{ and } p_z\in Y_{N-1}\}$ is a matching in $G$. 
Set $Y_{N-1}'=Y_{N-1}-V(M_N')$ and repeat the above process to vertices in $Y_{N-1}'$, we obtain matchings $M_{N-1}$ and $M_{N-1}'$.  Set $Y_{N-2}'=Y_{N-2}-V(M_{N-1}')$. 
This process is continued until we obtain $Y_1'=Y_1-V(M_2')$. We can still find a maximum matching $M_1$ in $G[Y_1']$ and a stable set $Z_1=Y_1'-V(M_1)$. 
By \cref{pq}, for every $u_i\in Y_0$,  there is at most one vertex in $Z_1$ satisfying $i=\min \{j : v\in Z_1 \text{ and } vu_j\in E(G) \}$, denote such a vertex  by $v_i$ (if any). By \cref{i+1}, $v_iu_{i+1}\in E(G)$. We retain the two edges $v_iu_i$ and $v_iu_{i+1}$ and delete all the other edges incident to $v_i$ for every possible $v_i$; let $L$ be the resulting subgraph of $G[Y_0\cup Z_1]$.

\begin{claim}\label{L}
	 $L$  is connected and claw-free. 
	
\end{claim}
\begin{proof}
Note that $V(L)=Y_0\cup Z_1$ as $Z_1\subseteq N(Y_0)$. By the construction of $L$, we have $L$ is connected and $d_L(v_i)=2$ for each $v_i\in Z_1$. Thus no vertex in $Z_1$ can be a claw center. 
For each vertex $u_i\in Y_0$, we have $d_L(u_i)\le 4$ as the possible neighbors of $u_i$ are $u_{i-1}, u_{i+1}\in Y_0$ and $v_{i-1}, v_i\in Z_1$. As
$N_L(v_{i-1})=\{u_{i-1},u_i\}$ and $N_L(v_i)= \{u_i,u_{i+1}\}$, $u_i$ cannot be a claw center either. 
\end{proof}

\begin{claim}\label{def(Y)}
$\kd(G[Y])\le 1$.	
\end{claim}
\begin{proof}
By Claim~\ref{L} and \cref{claw-free}, we have $\kd(L)\le 1$. Since all vertices in  $Y\setminus (Z_1\cup Y_0)$ has been saturated by the matching $\bigcup_{i=1}^N(M_i\cup M_i')$, where we set $M_1'=\emptyset$ for convenience. Therefore, we have $\kd(G[Y])\le 1$.  
\end{proof}

Since $G[X_0]$ consists of two disjoint paths, $G[X]$ is composed of at most two components. According to  Claim~\ref{CL: X_2n0+1}, each component has diameter at most $3n$. 
By \cref{finite level}, $\kd(G[X])$ is bounded by a function $b(n)=f(3n, n)$. Together with $\kd(G[Y])\le 1$, we know that $\kd(G)$ is bounded. The proof of \cref{main} is complete.

\end{proof}

\section{Remarks and Discussions}
According to the proofs of Theorems~\ref{main} and \ref{finite level}, every $\mathcal{H}$-free graph $G$ with $\mathcal{H}\in \B\kd$ has deficiency bounded by a function dependent on the Ramsey number $R(n, \alpha_i)$ and its diameter, where $\alpha_i$ is defined recursively dependent on $R(n,\alpha_{i-1})$. Seemingly, it is difficult to determine the exact value of $\kd(G)$ for an $\mathcal{H}$-free graph $G$ with $\mathcal{H}\in \B\kd$.
However, it will be very interesting to determine the family $\mathcal{H}\in d$-$\B\kd$ for a fixed integer $d$ as in Theorem~\ref{bone}. 

With a similar discussion as in Theorem~\ref{pair}, we have further determined $\mathcal{H}\in\B\kd$ with $|\mathcal{H}|=3$.
However, it is interminable to write down.

\subsection*{Acknowledgements}

\noindent
This work was supported by the National Key Research and Development Program of China (2023YFA1010200), the National Natural Science Foundation of China (No. 12071453),  and the Innovation Program for Quantum Science and Technology (2021ZD0302902).


\end{document}